\documentclass[11pt]{amsart}
\usepackage[latin9]{inputenc}
\usepackage{color}
\usepackage{cancel}
\usepackage{amsbsy}
\usepackage{amstext}
\usepackage{amsthm}
\usepackage{amssymb}
\usepackage{mathtools}
\PassOptionsToPackage{normalem}{ulem}
\usepackage{ulem}

\makeatletter
\numberwithin{equation}{section}
\numberwithin{figure}{section}

\usepackage{tikz-cd}
\usepackage{fancyhdr}
\usepackage{latexsym}
\usepackage{amscd}\usepackage{amsfonts}\usepackage{pstricks}
\usepackage{young}
\usepackage{enumitem}
\usepackage{amsthm}
\usepackage{mathrsfs}

\usepackage{ifthen}
\usepackage{longtable}
\usepackage{todonotes}
\usepackage{marginnote}

\renewcommand{\subsection}[1]{\vspace{3mm}\refstepcounter{subsection}\noindent{\bf \thesubsection. #1.} }

\renewcommand{\subsubsection}[1]{\vspace{3mm}\refstepcounter{subsubsection}\noindent{\bf \thesubsubsection. #1.} }

\numberwithin{equation}{section}

\newtheorem{theorem}{Theorem}
\newtheorem{lemma}[theorem]{Lemma}

\newtheorem*{theorem*}{Theorem}

\theoremstyle{definition}

\theoremstyle{remark*}
\newtheorem*{remark*}{Remark}

\newtheorem*{example*}{Example}

\def\min{\mathop{\mathrm{min}}}

\makeatother

\begin{document}
\title{Campana's orbifold conjecture for numerically equivalent divisors}
 
\author{Min Ru}
\address{
 Department  of Mathematics\newline
\indent University of Houston\newline
\indent Houston,  TX 77204, U.S.A.} 
\email{minru@math.uh.edu}
\author{Julie Tzu-Yueh Wang}
\address{Institute of Mathematics, Academia Sinica \newline
\indent No.\ 1, Sec.\ 4, Roosevelt Road\newline
\indent Taipei 10617, Taiwan}
\email{jwang@math.sinica.edu.tw}
\thanks{2020\ {\it Mathematics Subject Classification}: Primary 32H30; Secondary 32Q45 and 30D35}
\thanks{The research of Min Ru is supported in part by Simon Foundations grant award \#531604 and \#521604.
The  second-named author was supported in part by Taiwan's NSTC grant  113-2115-M-001-011-MY3.}

%\thanks{The   fourth-named author was supported in part by Taiwan's NSTC grant  110-2115-M-001-009-MY3.}

\begin{abstract} We prove  the following version of the Campana's orbifold conjecture: Let $X$ be a complex  non-singular projective variety of dimension $n$.  Let  $D_1,\hdots,D_{n+1}$ be  $\mathbb Z$-linearly independent effective divisors in ${\rm Div}(X)$ and $D:=D_1+\cdots+D_{n+1}$  be a normal crossing divisor of $X$.  Assume furthermore that they are numerically parallel.
Let $\Delta=\sum_{i=1}^{n+1} (1-m_i^{-1}) D_i$ and let  $f:\mathbb C\to (X,\Delta) $ be an orbifold entire curve. 
Then, there exists  a positive integer  $\ell$    such that,  the orbifold $ (X,\Delta_{\ell}) $ is of general type, where $\Delta_{\ell}=\sum_{i=1}^{n+1} (1-\frac1{\ell})D_i$, and  if $f$ has multiplicity at least $\ell$ along $D_i$, $1\le i\le n+1$,  then $f$ must be algebraically degenerate.
\newline
\begin{center}
{\it  Dedicated to Professor Josip Globevnik for his upcoming 80th's birthday in 2025} 
\end{center}\end{abstract}

\maketitle
\baselineskip=16truept

%%%%%%%%%%%%%%%%%%%%% ENDTOP MATTER %%%%%%%%%%%%%%%%%%
 \section{Introduction}\label{sec:intro}

\bigskip We start by recalling Campana's orbifold notion.  Let $X$ be a projective variety and $D$ be an effective Cartier divisor on $X$.
 Instead of considering the complement  $X\backslash  D$  in the  Brody's hyperbolicity case, one examines the  entire curve  $f:\mathbb C\to X$ that is ramified over $D$. According to Campana \cite{Camp}, 
	 an {\it orbifold divisor} is  defined as  $\Delta: =\sum_{Y\subset X} (1-m_{\Delta}^{-1}(Y))\cdot Y$, 
	 where $Y$ ranges over all irreducible divisors of $X$,   and $m_{\Delta}(Y ) \in [1, \infty)  \cap {\Bbb Q}$ with all but finitely many values equal to one.
	 An {\it orbifold entire curve} $f: {\Bbb C}\rightarrow  (X,\Delta)$, with $\Delta=\sum_i (1-m_i^{-1}) D_i$,  is a holomorphic  map such that $f({\Bbb C})\not\subset \mbox{Supp}(\Delta)$ and ${\rm mult}_{z_0}(f^*D_i) \ge  m_i$ for all $i$ and 
for every $z_0 \in {\Bbb C}$  with $f(z_0) \in  D_i$.   Let $\ell$ be a positive integer.  It is convenient to say that the multiplicity of the orbifold pair $ (X,\Delta)$  is at least $\ell$ if each  $ m_i\ge \ell$.  Note that we do not assume $D_i$ to be irreducible  or reduced.

As a consequence of the Main Theorem in \cite{RWorbifold2024},  the following result was obtained by the authors (see also \cite[Theorem 1.2]{GW22}).
\begin{theorem*}[Ru-Wang \cite{RWorbifold2024}, Theorem 1]\label{GG_conj}
		Let 
		$F_i$, $1\le i\le n+1$, be homogeneous irreducible polynomials of positive degrees  in $\mathbb{C}[x_0,\hdots,x_n]$.
		Assume that the hypersurfaces  $D_i:=[F_i=0]\subset\mathbb P^n(\mathbb C)$,  $1\le i\le n+1$,   intersect transversally.  Let $\Delta$ be an orbifold divisor of ${\Bbb P}^n(\mathbb C)$ with ${\rm Supp}(\Delta)=D_1+\cdots+D_{n+1}$.  Assume that  $\deg \Delta>n+1$.
		Then there exist  an effectively computable  positive integer  $\ell$ and a proper Zariski closed subset $W$ of  ${\mathbb P}^n(\mathbb C)$ such that if the orbifold $ ({\mathbb P}^n(\mathbb C),\Delta) $ has  multiplicities at least $\ell$, then the image of any   orbifold entire curve  $f: {\Bbb C}\to ({\mathbb P}^n(\mathbb C),\Delta) $ must be  contained in $W$.
	\end{theorem*}

The purpose of this paper is to consider,  instead of  ${\Bbb P}^n({\Bbb C})$,  a general  complex non-singular projective variety $X$ and  where the given divisors are  numerically parallel on $X$.

Let $X$ be a complex (non-singular) projective variety of dimension $n$.  A set of  effective divisors  $D_1,\hdots,D_l$, $l\ge 2$, on $X$
is said to be {\it numerically parallel} if there exist positive integers $d_i$ such that $d_jD_i\equiv d_iD_j$, where $``\equiv"$ means that two divisors are numerically equivalent on $X$.
Note that on ${\Bbb P}^n({\Bbb C})$, all effective divisors are numerically parallel  (indeed they are always linearly parallel, i.e.  $d_jD_i\sim  d_iD_j$, where ``$\sim$" means that two divisors are linearly equivalent on $X$).

Denote by ${\rm Div}(X)$ the group of effective divisors on $X$.
We say that {\it $D_i, 1\leq i\leq l$, are  $($or the family $\{D_i\}_{i=1}^l$ is$)$  $\mathbb Z$-linearly independent divisors in ${\rm Div}(X)$} if, for every $1\leq i, j\leq l$,
$$\mbox{Supp} D_j\not\subset \bigcup_{i\not=j} \mbox{Supp} D_i.$$

The main theorem of this paper is as follows.
\begin{theorem}\label{MainTheorem}
Let $X$ be a complex  non-singular  projective variety of dimension $n$.  Let  $D_1,\hdots,D_{n+1}$ be  $\mathbb Z$-linearly independent effective divisors in ${\rm Div}(X)$ such that $D:=D_1+\cdots+D_{n+1}$  is of normal crossing  and $(X, D)$ is of log general type. 
 Assume furthermore that they are numerically parallel. Let $\Delta$ be an orbifold divisor with ${\rm Supp} \Delta={\rm Supp} D$. Let $f:\mathbb C\to (X,\Delta) $ be an orbifold entire curve. 
Then there exists  an integer $\ell>0$   such that,  the orbifold $ (X,\Delta_{\ell}) $ is of general type, where $\Delta_{\ell}=\sum_{i=1}^{n+1} (1-\frac1{\ell})D_i$, and  if $f$ has multiplicity at least $\ell$ along $D_i$, $1\le i\le n+1$,  then $f$ must be algebraically degenerate, i.e. the image of $f$ is contained in a proper subvariety of $X$.
\end{theorem}  
Let $D:={\rm Supp}(\Delta)= D_1+\cdots+D_{l}$.
		Recall that the pair $(X, D)$ is said to be   of {\it log-general type} if the divisor $K_X+D$ is big.
		It's clear that  $(X, D)$ is of log-general type if  $(X,\Delta)$ is   of  general type, i.e. $K_X + \Delta$ is big.
		On the other hand, the condition that $(X, D)$ is of log-general type implies that $(X,\Delta)$ is   of  general type if $ (X,\Delta) $ has a sufficiently large multiplicity along  each $D_i$, $1\leq i\leq l$.  $($See \cite[Corollary 2.2.24]{LazarsfeldI}.$)$

\section{The orbifold version of the  logarithmic Bloch-Ochiai Theorem}

In this section, we will establish the orbifold version of  the  logarithmic Bloch-Ochiai Theorem and an orbifold analogue of Laurent's fundamental result \cite{Laurent}, which describes the Zariski closure of subsets of finitely generated subgroups of $\mathbb G_m^n(\overline {\mathbb Q})$.

 We use the standard notations in Nevanlinna theory (see \cite{ru2021nevanlinna} or  \cite{noguchibook}).   Let $X$ be  a complex projective variety and 
 $f: {\Bbb C}\rightarrow X$ be a holomorphic map.  We fix a  ample line bundle $L$ on $X$, and let $T_{f, L}(r)$ be the characteristic function of $f$ with 
 respect to $L$. Most of time in this paper, we just write $T_f(r)$ by dropping the line bundle $L$. We also use $S_f(r)$ to denote the term which is 
 $\leq_ {\rm exc}  O(\log T_f(r))$ where ``$\leq_ {\rm exc} $" means that the inequality ``$\leq$" holds for all $r\in (0, +\infty)$ except for a set $E$ with finite measure. 
 
  \noindent \begin{lemma}[Lemma 4.7.1 in  \cite{noguchibook}]\label{logjet} Let $X$ be  a complex projective variety and let $D$ be a reduced divisor on $X$. Let $f:  {\Bbb C}\rightarrow X$ be a holomorphic map such that $f({\Bbb C})\not\subset D$. 
  Let $\omega$ be a  log $k$-th jet differential along $D$, and write $\xi(z):=\omega(j_k(f))(z)$.
  Then 
  $$m(r, \xi) =S_f(r),$$
  where $m(r, \xi)={1\over 2\pi}\int_0^{2\pi}  \log^+|\xi(re^{i\theta})| d\theta$.    
\end{lemma}

We need some preparations.  Let $A$  be a semi-abelian variety of dimension $n$ over ${\Bbb C}$; that is, $A$ is a
complex algebraic group of dimension $n$ which is an extension of an Abelian variety $A_0$ by  an algebraic torus $({\Bbb C}^*)^t$ ($t\in {\Bbb Z}_{\ge 0}$),
$$0 \to ({\Bbb C}^*)^t\to A\to A_0\to 0.$$
We have a decomposition 
$$J_k(A)=A\times {\Bbb C}^{nk},$$
where $J_k(A)$ is the the $k$-jet space of $A$.
Let $X\subset A$ be an irreducible algebraic subvariety, we have the stabilizer group 
(with respect to the natural $A$-action) 
$$\mbox{St}_A(X)=\{a\in A|~a+x\in X~ \mbox{for  all}~x\in X\}$$ 
and the identity component $\mbox{St}_A(X)^0$. The $\mbox{St}_A(X)^0$ is a connected closed algebraic group of $A$.
Let 
\begin{equation}\label{rho}
\rho_k: J_k(A)= A\times {\Bbb C}^{nk}\rightarrow {\Bbb C}^{nk} 
\end{equation}
be the natural projection.

\begin{lemma}\label{Lie}  Let $A$  be a semi-abelian variety and $X\subset A$ be an irreducible algebraic subvariety
Let $\rho_k$ be the map defined as above.  Assume that  $\mbox{\rm St}_A(X)^0=\{0\}$. Then there exists  a sufficiently large  $k_0>0$ such that for every
  $k\ge k_0$, the differential $d\rho_k|_{J_k(X)}$  is   generically injective. 
\end{lemma} 
Note that  this lemma is indeed a special case of Lemma 6.2.4 in  \cite{noguchibook}.  

\begin{proof}  The proof basically just repeats the argument of the proof of Lemma 1.2  in \cite{No08} or Lemma A9.1.7 in \cite{ru2021nevanlinna}.

We fix arbitrarily a point $y_0\in X$. By  Lemma 6.2.2 in  \cite{noguchibook}, there is an algebraically non-degenerate holomorphic map $f: \bigtriangleup(1)\to X$ with $f(0)=y_0$. Let $J_k(f):  \bigtriangleup(1)\rightarrow    J_k(X)$ denote the k-th jet lifting of $f$, where $J_k(X)$ is the the $k$-jet space of $X$.
 Set $y_k = J_k(f)(0) \in  J_k(X)$. Let $x_0\in A$ be the base point of $y_k$ by the projection to the base space $J_k(X)\rightarrow A$.
 First note that 
 $${\bf T}_{y_k}(J_k(X))\subset {\bf T}_{y_k}(A\times {\Bbb C}^{nk})\cong {\bf T}_{x_0}(A)\oplus {\Bbb C}^{nk},$$
 where ${\bf T}(*)$ stands for the holomorphic tangent space. 
 Because of the definition of $\rho_k|_{J_k(X)}$, we have 
\begin{equation}\label{N1}
 \mbox{Ker}(d\rho_k|_{J_k(X)}(y_k)) \subset   {\bf T}_{x_0}(A) \oplus  O' \cong {\bf T}_{x_0} (A).
 \end{equation}
 Note that,  for  $k'\ge k$, 
$$ \mbox{Ker}(d\rho_{k'}|_{J_{k'}(X)}(y_{k'} )) \subset   \mbox{Ker}(d\rho_k|_{J_k(X)}(y_k)),$$
thus if $ \cap_{k\ge 1} \mbox{Ker}(d\rho_k|_{J_k(X)}(y_k))=\{0\}$ then $\mbox{Ker}(d\rho_{k_0}|_{J_{k_0}(X)}(y_{k_0}))=\{0\}$ for some $k_0>0$. This will imply that the differential $d\rho_{k_0}|_{J_{k_0}(X)}$ is injective at $y_{k_0}$.
Therefore it remains to prove that 
\begin{equation}\label{W}\bigcap_{k\ge 1} \mbox{Ker}(d\rho_k|_{J_k(X)}(y_k))=\{0\}.\end{equation}
Suppose that 
$\cap_{k\ge 1} \mbox{Ker}(d\rho_{k}|_{J_{k}(X(f))}(y_{k}))\not =\{0\}$. Take a non-trivial vector $$v\in  \cap_{k\ge 1}\mbox{Ker}(d\rho_{k}|_{J_{k}(X(f))}(y_{k})).$$
Let $h\in {\mathcal I}(X)_{y_0}$ be a germ of holomorphic functions in the ideal sheaf of $X$ at $y_0$. Thus $v$, considered as a vector field on $A$, satisfies
$${d^k\over dz^k}\big |_{z=0}vh(f(z))=0, k=0, 1, 2, ...$$
Therefore $z\mapsto vh(f(z))$ has zero power series development, i.e. 
$vh(f(z))\equiv 0$
in a neighborhood of $0$. Since $h$ was chosen arbitrarily, we obtain $vh(f(z))\equiv 0$ near $0$ for all $h\in {\mathcal I}(X)_{y_0}$. Hence 
$v$ is tangent to $X$ at every  point of $f(\bigtriangleup(1))$. Since $f$ is Zariski dense in $X$, $v$ is everywhere tangent to $X$, so $X$ is invariant by the action of the one-parameter subgroup generated by $v$.
This  contradicts with the fact that $\mbox{St}_A(X)^0=\{0\}$.  Then (\ref{W}) holds. This finishes the proof.
\end{proof}

  Let $\bar{A}$ be the compactification of $A$
 such that $\partial A = \bar{A} \backslash  A$, which is a divisor on $\bar{A}$ with only simple normal crossings.  Thus
 we have the logarithmic $k$-th  jet bundle $J_k(\bar{A}; \log\partial A)$ over $\bar{A}$ along $\partial A$, and 
  we can trivialize it through the differential forms $dz_1, \dots, dz_n$
on $A$ where $(z_1, \dots, z_n)$ is the coordinates of ${\Bbb C}^n$. 
Then we get a decomposition $$J_k(\bar{A}; \log\partial A)={\bar A}\times {\Bbb C}^{nk}.$$ 
Let
$${\tilde \rho_k}:  J_k(\bar{A}; \log\partial A)={\bar A}\times {\Bbb C}^{nk} \rightarrow {\Bbb C}^{nk}.$$
Note that ${\tilde \rho_k}|_{J_k(X)}=\rho_k|_{J_k(X)}$ where $\rho_k$ is defined in (\ref{rho}),  so we still write ${\tilde \rho_k}$ as $\rho_k$.

\begin{theorem}\label{Nog-Ru} Let $A$ be a semi-abelian variety and let $X\subset A$ be an algebraic variety with $\mbox{\rm St}_A(X)^0=\{0\}$. Let $\bar{A}$ be the compactification of $A$ and $\bar{X}$ be the compactification in $\bar{A}$.
Let $f: {\Bbb C}\to {\bar X}$ be an algebraically non-degenerate map. Assume that $D={\bar X}\backslash X$ $($i.e. $D:=\tau^*(\partial A)|_{\bar{X}}$, where $\tau: A\to \bar{A}$ is an inclusion map$)$. Then there exists a constant  $\kappa>0$, such that 
$$\kappa T_f(r) \leq N_f^{(1)}(r, D)+S_f(r).$$
\end{theorem}

\begin{proof} By the lemma above, since $ {\rm St}_A(X)^0=\{0\}$, there exists  a sufficiently large  $k_0>0$ such that for every 
  $k\ge k_0$, the differential $d\rho_k|_{J_k(X)}$    is generically injective. We fix such  a $k$.
  
   Let $\{\phi_l\}_{l=1}^s$ be a transcendental basis of the 
rational function field  ${\Bbb C}({\bar X})$ of ${\bar X}$.  
Then by Theorem 2.5.18 in  \cite{noguchibook} there is some constant $C>0$ such that 
\begin{equation}\label{chara}
C^{-1}T_f (r) \leq T_{f}(r, \{\phi_j\})\leq  CT_{f}(r)+O(1),
\end{equation}
where $T_{f}(r, \{\phi_l\})=\sum_{l=1}^s T_{\phi_l(f)}(r).$
 
 Let $\omega_1, \dots, \omega_q$ be the bases of $H^0({\bar X}, \Omega_{\bar{X}}^1(\log D))\subseteq H^0({\bar A}, \Omega_{{\bar A}}^1(\log \partial A)).$
Put $$\zeta_i^j(z)=d^j\omega_i(J_k(f)(z)).$$ From  logarithmic derivative  Lemma (Lemma \ref{logjet}),
$$m(r, \xi_i^j)\leq S_f(r),$$
and thus, by the First Main Theorem, 
\begin{equation}\label{lo}
T(r, \xi_i^j )\leq (j+1)N_f^{(1)}(r, D) +S_f(r).
\end{equation}
 Let $Z$ be the Zariski closure of $J_k(f)$ in $J_k({\bar X})$ and let
 $V$ be the Zariski closure of $\rho_k|_{J_k(f)({\Bbb C})}$ in ${\Bbb C}^{nk}$.  Since  $d\rho_k$ has maximal rank at $J_k(f)(0)$, it follows that $\dim Z=\dim V$ and that 
  $\rho_k|_Z: Z\rightarrow V$ is   dominant. We can regard  $\phi_l$ as rational functions on $Z$ through $p|_Z: Z \rightarrow {\bar X}$. 
  By Theorem A9.1.8  in \cite{ru2021nevanlinna},  there are 
algebraically relations, for each $l$, 
$$P_{0}(d^j\omega_i)\phi_l^{ d_l}+  \dots + P_{d_l}(d^j\omega_i)\equiv 0,$$
$$\mbox{with}~~~P_0(d^j\omega_i)\not\equiv 0.$$
Therefore $$P_{0}(\xi_i^j)f^*\phi_l^{d_l}+  \dots + P_{d_l}(\xi_i^j)\equiv 0,$$
$$\mbox{with}~~~P_0(\xi_i^j)\not\equiv 0.$$
This implies, by Lemma 2.5.15 in     \cite{noguchibook} and  (\ref{lo})  that 
$$T_f(r,  \{\phi_l\})  \leq  c_1 N_f^{(1)}(r, D)+S_{f}(r) $$
for some positive constant $c_1$.
This together with (\ref{chara}) implies that 
$$T_f(r)   \leq Cc_1 N_f^{(1)}(r, D)+S_f(r).$$
This proves the theorem.
  \end{proof}
\begin{theorem}[Orbifold  logarithmic Bloch-Ochiai Theorem, see Theorem 4.8.17 in  \cite{noguchibook}]\label{logB} Let $X$ be a   non-singular projective  variety and let $D$ be an effective  reduced divisor on $X$.   
Let $\Delta$ be an orbidold divisor on $X$ with $ {\rm Supp} \Delta =  D$.
Let $f: {\Bbb C}\rightarrow (X,\Delta)$ be an orbifold  holomorphic curve.  Assume that $\bar{q}(X\backslash D)>\dim X$. Then there exists a positive integer $\ell$ such that if the multiplicity of $\Delta$ is at least $\ell$, then
  $f$ is algebraically degenerate. 
\end{theorem}
\begin{proof}  Assume  $f$ is algebraic non-degenerate. 
 Let $$\alpha: X \backslash D \rightarrow A_{X\backslash D}$$ be the quasi-Albanese map.
This extends to a holomorphic map $\bar{\alpha}: X \rightarrow \bar{A}_{X\backslash D}$, where $\bar{A}_{X\backslash D}$ is the compatification of 
$A_{X\backslash D}$. Let $g:=\bar{\alpha}(f): {\Bbb C}\rightarrow  \bar{A}_{X\backslash D}$. Denote by ${\bar Y}$ the Zariski closure of $g({\Bbb C})$ in $\bar{A}_{X\backslash D}$,  
and let  $Y=\bar{Y}\cap A_{X\backslash D}$.  Then $\bar{Y}\setminus Y\subset \bar{A}_{X\backslash D}\setminus A_{X\backslash D}$.
 By assumption, $\dim Y =\dim X<\bar{q}(X\backslash D)=\dim A_{X\backslash D}$.   Since   $ A_{X\backslash D}$ is generated by the image of $\alpha(X \backslash D)$, it implies that  $Y$ is not contained in any translate of subgroup of $A_{X\backslash D}$.

 Set $A_1=A_{X\backslash D}/{\rm St}_{A_{X\backslash D}}(Y)^0$.  Then $\dim A_1>0$ and we know that $A_1$ is a semi-abelian by the arguments in the proof of Theorem 4.8.17 in  \cite{noguchibook} (alternatively, use Proposition 5.1.26 in \cite{noguchibook}).
 Let $$\beta: A_{X\backslash D} \rightarrow A_1$$ be the quotient map and let $Y_1= {\beta}(Y)$.  Then ${\rm St}_{A_1} (Y_1)^0=\{0\}$ in $A_1$  and 
$\beta$ extends to a meromorphic map $\bar\beta:\bar A_{X\backslash D}\to \bar A_1$ between compactifications such that $\bar{\beta}(\bar Y)={\bar Y_1}$.  Thus $\bar{Y_1}\setminus Y_1\subset  \bar{A_1}\setminus A_1$. 
We apply Theorem \ref{Nog-Ru} to $h=\bar{\beta}(g): {\Bbb C}\rightarrow \bar{Y_1}$ (with $Y_1\subset A_1$) to get, for some constant $\kappa>0$, 
\begin{align}\label{char}
\kappa T_h(r) \le_{\rm exc}   N^{(1)}_h(r, E),
\end{align}
where $E=\bar{Y_1}\setminus Y_1$.

Let $\Delta_1$ be the reduced divisor of $E:=\bar{Y_1}\setminus Y_1$.
For the given   orbifold  holomorphic curve $f: {\Bbb C}\rightarrow (X,\Delta)$   with multiplicity bigger than a sufficiently large integer $\ell$ to be determined later, 
first note that  $h=\bar{\beta}\circ \bar{\alpha}\circ f:\mathbb C\to \bar Y_1$.  Furthermore, for any irreducible component   $W$  of  $E$,  note that 
$ \bar{\alpha}^*(\bar{\beta}^*W)$ is contained in the support of $D$, we have,
 for every $z_0 \in {\Bbb C}$  with $h(z_0) \in  W$, 
 $ {\rm mult}_{z_0}(h^*W)= {\rm mult}_{z_0}(f^* \bar{\alpha}^*\bar{\beta}^*W)\ge  \ell $. Hence,
\begin{align}\label{mutiW}
N^{(1)}_h(r, E)\le  \frac1\ell T_{h, E}(r).
\end{align}
By combing the above inequality with (\ref{char}), we get 
$$\kappa T_h(r) \le_{\rm exc}   \frac1\ell T_{h, E}(r).$$
Note that  $T_h(r)\ge c  T_{E, h}(r)$ for some positive constant $c$ since  $T_h(r)=T_{h, L'}(f)$ for some ample divisor $L'$ on $\bar{Y_1}$.
We get a contradiction by allowing $\ell$ sufficiently larger than $(c \cdot\kappa)^{-1}$.
\end{proof}

By using the theorem above, we can get the following variant of Theorem  4.9.7  in \cite{noguchibook}, whose proof is exactly the same as the  proof of Theorem 4.9.7  in \cite{noguchibook}. 
Let \( \mathrm{NS}(X)\) denote the  N\'eron-Severi group, i.e.,
$$\mathrm{NS}(X) =\mathrm{Pic}(X)/\mathrm{Pic}^0(X).$$ Let $D_i, 1\leq i\leq l,$ be distinct  divisors on $X$.  Denote by 
 \( \mathrm{rk} \{ D_i \}_{i=1}^{l} \)  the rank of the subgroup  $\sum_{i=1}^l  {\Bbb Z}\cdot c_1(D_i)$ generated by $c_1(D_i)\in H^2(X, {\Bbb Z})$, $1\leq i\leq l$, in  $H^2(X, {\Bbb Z})$.
 Note that 
 $$\mathrm{rk} \{ D_i \}_{i=1}^{l} \leq \mathrm{rank}_{\Bbb Z} \mathrm{NS}(X).$$
\begin{theorem}[c.f.   Theorem 4.9.7  in \cite{noguchibook}]\label{newtheorem} Denote by $q(X):= h^1(X, \mathcal{O}_X) $. Let $D_i, 1\leq i\leq l$, be distinct  reduced divisors on $X$. 
 Assume that  $l>\dim X+\mathrm{rk} \{ D_i \}_{i=1}^{l} -q(X)$.  Let $\Delta$ be an orbifold divisor on $X$ with $ {\rm Supp} \Delta =  D$, where $D=D_1+\cdots+D_l$.
Let $f: {\Bbb C}\rightarrow (X,\Delta)$ be an orbifold holomorphic curve. Then there exists a positive integer $\ell$ such that if the multiplicity of $\Delta$ is at least $\ell$, then
  $f$ is algebraically degenerate. 
 \end{theorem}

\begin{proof}
By  Lemma 4.9.6  in \cite{noguchibook}, we have $\bar{q}(X\backslash D)\ge q(X)+l-\mathrm{rk} \{ D_i \}_{i=1}^{l}.$  If  $l>\dim X+\mathrm{rk} \{ D_i \}_{i=1}^{l} -q(X)$, then $\bar{q}(X\backslash D)>\dim X$.  Consequently, the assertion follows from Theorem \ref{logB}.
\end{proof}

We recall the following theorem from \cite{noguchibook}.  We will use this theorem for the case that $A=\mathbb G_m^n$ and $\bar A=\mathbb P^n({\Bbb C})$.

\begin{theorem}[Theorem 5.3.23 in \cite{noguchibook}]\label{newtheorem2}  Let $\bar{A}$ be a semi-toric variety and let $Z$ be a closed irreducible analytic subset of $A$ such that its closure $\bar{Z}$ in ${\bar A}$ is an analytic subset. Then $Z$ is a translate of a subgroup iff $\bar{q}(Z)=\dim Z$.
\end{theorem}
\begin{proof} We recall that $\bar{q}(Z)$ equal to the dimension of the quasi-Albanese variety of $Z$. If $Z$ is a translate of a subgroup, then it is its own quasi-Albanese variety, implying that $\bar{Z}=\dim Z$. If 
$Z$ is not a translate of a subgroup we regard the subgroup $B$ of $A$ generated by all $p-q (p, q\in Z)$. This is a quasi-algebraic subgroup (see  \cite{noguchibook}) due to Proposition 5.3.17 in \cite{noguchibook} and in particular a semi-torus. Hence we obtain a surjection from the 
quasi-Albanese variety of $Z$ onto $B$, implying $\bar{q}(Z)\ge \dim B>\dim Z$. 
\end{proof}

 We will need the following theorem, which is  a consequence of  the orbifold logarithmic   Bloch-Ochiai Theorem (Theorem \ref{logB}).

\begin{theorem}\label{Ru}   Let $\Delta$ be an orbifold divisor on $\mathbb P^n({\Bbb C})$ with {\rm Supp}$\Delta =H$, where $H$ is the union of the coordinate hyperplanes.
Let $g : {\Bbb C}\rightarrow (\mathbb P^n({\Bbb C}),\Delta)$ be an orbifold entire curve.  
Then there exists a positive integer $\ell$ such that if the multiplicity of $\Delta$ is at least $\ell$, then   the Zariski closure of $g({\Bbb C}\backslash g^{-1}(H))$ in $\mathbb G_m^n$  is a translate of a subgroup of $\mathbb G_m^{ n }: = {\Bbb P}^n({\Bbb C})\backslash H$. 
\end{theorem}

\begin{remark*}
This theorem  in the unit case follows from the logarithmic   Bloch-Ochiai Theorem(Theorem 4.8.17 in  \cite{noguchibook}) and Theorem \ref{newtheorem2}, under the assumption that  the image of $g$ is contained  in $\mathbb G_m^n$.  Moreover,  the unit case  serves as the complex analogue of  Laurent's fundamental result \cite{Laurent}, which describes the Zariski closure of subsets of finitely generated subgroups of $\mathbb G_m^n(\overline {\mathbb Q})$.
\end{remark*}

 \begin{proof}   
Let $g=(g_0,g_1,\hdots,g_n)$, where $g_0,\hdots,g_n$ are entire functions with no common zero.
We first treat the case when each $g_i$ is a unit, meaning it has neither zeros nor poles.
By normalizing with respect to $g_0$, we regard  $g$ as a map from $\mathbb C$ to $ \mathbb G_m^{n}$.  
Let $Z$ be the Zariski closure of $g(\mathbb C)$ in $\mathbb G_m^{n}$.  For our purpose, we may assume that $Z$ is neither a point nor the full $\mathbb G_m^{n}$.  If $Z$ is not a translate of subgroup of $G_m^{n}$, then $\bar q(Z)>\dim Z$ by Theorem \ref{newtheorem2}.  By the Bloch-Ochiai Theorem (see  \cite[Corollary 4.8.18]{noguchibook}), this forces $g:\mathbb C\to Z$ to be algebraically degenerate, contradicting to our choice of $Z$.
 
We now consider the general case.  Let $Z$ be the Zariski closure of $g({\Bbb C}\backslash  g^{-1}(H))$ in  $ \mathbb G_m^{n}$.  
 We note that $Z$ is not empty since  $g(\mathbb C)\not\subset H$.   Let $\bar Z$ be the Zariski closure of $Z$ in $\mathbb P^n({\Bbb C})$.
Write  $\Delta=\sum_{i=0}^n(1-\frac1 {m_i}) H_i$, where $H_i=[x_i=0]$.  
We first consider the case where  $H_i\cap \bar Z=\emptyset$  for all $0\le i\le n$.  In this case, each $g_i$ must be a unit, and this situation has already been treated at the beginning of the proof. Now,  assume that  at least one intersection $H_i\cap \bar Z$ is nonempty for some $0\le i\le n$. Since  $\bar Z$ is not contained in any $H_j$, it follows that   $H_j\cap \bar Z$ is a divisor of  $\bar Z$ whenever  $H_j\cap \bar Z$ is nonempty.
 Let $W_1,\hdots,W_k$ be  the irreducible components  appearing in some $H_j\cap \bar Z$ for $0\leq j\leq n$ with   the intersection $H_j\cap \bar Z$ being non-empty.
Define  
$$
\Delta_Z:=\sum_{j=1}^k (1-\frac 1{\min\{m_0,\hdots,m_n\}})W_j.
$$
For any $z_0 \in {\Bbb C}$ such that  $g(z_0) \in  W_j$, where $W_j$ is a component of $H_{j_0}\cap \bar Z$ for some $0\le j_0\le n$, we find
${\rm mult}_{z_0}(g^*W_j) \ge  m_{j_0}\ge \min\{m_0,\hdots,m_n\}.$
Therefore, $g: {\Bbb C}\rightarrow (\mathbb P^n({\Bbb C}),\Delta)$ induced an   orbifold curve $\tilde g =g:\mathbb C\to (\bar Z, \Delta_Z)$. 
 If $Z$ is not a translate of subgroup of $G_m^{n}$, then $\bar q(Z)>\dim Z$ by Theorem \ref{newtheorem2}.  We note that  $\bar q(Z)=\bar q(\bar Z\setminus \bar Z\cap H)$. 
 By Theorem \ref{logB} we find  a sufficiently large $\ell$ such that if the multiplicity of $ \Delta_Z$, i.e. $\min\{m_0,\hdots,m_n\}$,  is at least $\ell$, then $\tilde g$  is algebraically degenerate contradicting to the assumption of $Z$. 
\end{proof}
 
 \section{Proof of Theorem \ref{MainTheorem}}
 The following lemma is the orbifold version of \cite[Lemma 3.1]{LevinHuang}.  Its proof follows from their arguments by replacing  the use of a result of   Vojta~\cite{Vojtasemi} on integral points on subvarieties of semiabelian varieties by Theorem \ref{newtheorem}.

\begin{lemma}\label{LevinHuang}
Let $X$ be a complex non-singular projective variety of dimension $n$, and let $D_1, \ldots, D_{n+1}$ be $\mathbb{Z}$-linearly independent  effective divisors in $\mathrm{Div}(X)$. Assume additionally that these divisors are in general position and numerically parallel, i.e. there are integers $d_i, d_j$ such that $d_j D_i \equiv d_i D_j$ (numerical  equivalence) for all $i, j$.  Let $D:=D_1+\cdots+ D_{n+1}$.
Let $\Delta$ be an orbifold divisor on $X$ with $ {\rm Supp} \Delta = {\rm Supp} D$.  Suppose that there exists an algebraically non-degenerate
orbifold curve $f: {\Bbb C}\rightarrow (X,\Delta)$ with sufficiently large multiplicity  along  ${\rm Supp} \Delta$.
 Then 
 we have $q:= h^1(X, \mathcal{O}_X) = 0$, and  $d_j D_i \sim d_i D_j$ (linear equivalence) for all $i, j$.
Moreover, if $\mathrm{div}(\phi_i) = d_{n+1} D_i - d_i D_{n+1}$ for $i = 1, \ldots, n$, then the morphism
\[
\phi : X \setminus D \to \mathbb{G}_m^n, \quad \phi(P) = (\phi_1(P), \ldots, \phi_n(P)),
\]
is dominant.
\end{lemma}

\noindent{\bf Remark}.  The fact that, under the assumptions in the lemma above,  
$D_1, \dots, D_{n+1}$ being numerical equivalent implies that  they are indeed linear equivalent, seems to have independent interests. 

%\bigskip\noindent{\it Proof of Lemma \ref{LevinHuang}}
\begin{proof}[Proof of Lemma \ref{LevinHuang}] Assume that  $f: {\Bbb C}\rightarrow (X,\Delta)$ is algebraic non-degenerate orbifold curve, where the multiplicity of the orbifold  $(X,\Delta)$ is at least $\ell$ and  $\ell$ is an integer to be determined below.
By Theorem \ref{newtheorem}, if 
\begin{align}\label{dimA}
l> n+ \mathrm{rk} \{ D_i \}_{i=1}^{n+1}   -q,
\end{align} 
where $l$ is the number of irreducible component of $D$ and $\dim X=n$,  then  there exists a positive integer $\ell$ such that if the multiplicity of $\Delta$ is at least $\ell$,  $f$ must be degenerate.
Since  $D_1,\hdots,D_{n+1}$ are  $\mathbb Z$-linearly independent effective divisors   in ${\rm Div}(X)$, we have $l\ge n+1$.  Moreover,  under the assumption that the divisors $D_1, . . . , D_{n+1}$ are numerically parallel, it follows that  $\mathrm{rk} \{ D_i \}_{i=1}^{n+1}=1$. 
Hence, if $q>0$, then the inequality (\ref{dimA}) holds, forcing $f$ to be degenerate, which contradicts our assumption. Thus, we conclude that \( q = 0 \). In this case, the Albanese variety \( \mathrm{Alb}(X) \) is trivial. Since the Picard variety and Albanese variety are dual, the Picard variety \( \mathrm{Pic}^0(X) \) is also trivial. From the exact sequence
\[
0 \to \mathrm{Pic}^0(X) \to \mathrm{Pic}(X) \to \mathrm{NS}(X) \to 0,
\]
we have \( \mathrm{Pic}(X) \cong \mathrm{NS}(X) \). Since by assumption \( d_j D_i \equiv d_i D_j \) for all \( i \) and \( j \), it follows that \( d_j D_i \sim d_i D_j \) for all \( i \) and \( j \), i.e. $d_jD_i$ and $d_iD_j$ are linearly equivalent on $X$.

We now prove the morphism map
\[
\phi : X \setminus D \to \mathbb{G}_m^n, \quad \phi(P) = (\phi_1(P), \ldots, \phi_n(P)),
\]
is dominant, where $\mathrm{div}(\phi_i) = d_{n+1} D_i - d_i D_{n+1}$ for $i = 1, \ldots, n$.
By replacing $D_i$ with $(d/d_i)D_i$ where $d=\mbox{lcm}\{d_1, \dots, d_{n+1}\}$, we can assume, without loss of generality, that $d_i=1$ for $i=1, \dots, n+1$.
Take ${\Bbb P}^n({\Bbb C})$ as the compactification of $\mathbb{G}_m^n$.
Then the map $\phi$ extends to a  morphism 
\[
\bar \phi(P) : X\to \mathbb P^n({\Bbb C}). 
\]
Denote by $H=H_0+\cdots+H_n$ be the sum of coordinate hyperplanes $H_i=[x_i=0]$ in $\mathbb P^n({\Bbb C})$. Then
$
 \phi^{*}(H_i)=D_{i+1},\text{ for }  0\le i\le n+1$.
 Write $\Delta=\sum_{i=1}^{n+1} (1-m_i^{-1}) D_i$.   
Suppose that $(\bar\phi\circ f) (z_0)\in H_i$ for $z_0 \in {\Bbb C}$.
Then  
$$
{\rm mult}_{z_0}((\bar\phi\circ f)^*H_i) ={\rm mult}_{z_0}(  f^*(\bar\phi^*H_i))={\rm mult}_{z_0}(  f^*(D_{i+1})).
$$
Let  $\Delta_{H}:=(1-m_{n+1}^{-1}) H_0+\sum_{i=1}^{n} (1-m_i^{-1}) H_i$.
Then $\bar\phi\circ f:\mathbb C\to (\mathbb P^n({\Bbb C}),\Delta_H) $ is an entire orbifold curve.  Furthermore, we note that if   the multiplicity of $\Delta$ is at least $\ell$, then the multiplicity of $\Delta_H$ is also at least $\ell$.

Let ${\mathcal Y_{f}}:={\Bbb C}\backslash  f^{-1}( {\rm Supp}D)$.   
Suppose that  \( \phi \) is not dominant, then \( \phi( f({\mathcal Y_f})) \) is not Zariski dense in \( \mathbb{G}_m^n \). 
Let $Z$ be the Zariski closure of  \( \phi( f({\mathcal Y_f})) \)   in \( \mathbb{G}_m^n \).  Then by Theorem \ref{Ru}, there exists a positive integer $\ell_0$ such that if the multiplicity of $\Delta_H$ is at least $\ell_0$, then  $Z$ is a translate of a subgroup of $\mathbb G_m^{n }: = {\Bbb P}^n({\Bbb C})\backslash H$. 
This implies that  \( \phi(f({\mathcal Y_{f}})) \) is contained in a translate subgroup of $\mathbb G_m^n$.
 It follows that there is a nontrivial multiplicative relation \( \phi_1^{e_1} \ldots \phi_n^{e_n} = c \) with \( c \in \mathbb C^* \). Since \( \mathrm{div}(\phi_i) = d_{n+1} D_i - d_i D_{n+1} \), this yields a nontrivial linear dependence relation for \( D_1, \ldots, D_{n+1} \), contradicting our general setup. It follows that {\( \phi \)} is dominant.
\end{proof}

Recall the following theorem from  \cite{RWorbifold2024}.
\begin{theorem}[Ru-Wang  \cite{RWorbifold2024}] \label{finitemorphismtoric}
Let $Y$ be a nonsingular projective toric  variety, and let  $D_0 =Y\setminus \mathbb G_m^n$.
Let  $X$ be a nonsingular  complex   projective  variety with a finite  morphism $ \pi:  X\to Y$.  
Let $H:=\pi^*(D_0)$ and $R\subset X$ be the ramification divisor of $\pi$ omitting components from the support of $H$.
Assume that  $\pi(R)$ and $D_0$ are in general position on $Y$. 		
Then there exists a positive integer $\ell$ such that if the orbifold pair $(X,\Delta)$ is   of  general type with ${\rm Supp}(\Delta)=  {\rm Supp}(H)$ and  multiplicity at least $\ell$ for each component, then there exists a proper Zariski closed subset $W$ of $X$ such that  the image of any non-constant orbifold entire curve $f: {\Bbb C} \rightarrow (X,\Delta)$  must be contained in $W$. 
	\end{theorem}
 
%\noindent{\it  Proof of Theorem \ref{MainTheorem}.}
 \begin{proof}[Proof of Theorem \ref{MainTheorem}] We proceed by contradiction. Suppose that $f$ is algebraically non-degenerate.  
 For a positive integer $m$, we let $\Delta_m=\sum_{i=1}^{n+1} (1-\frac 1m)D_i.$
Since $D=D_1+\cdots+D_{n+1}$ is a normal crossing divisor, Lemma \ref{LevinHuang} implies that
there is a positive integer $\ell_1$ such that  if  $f$ has multiplicity at least $\ell_1$ along $D_i$, $1\le i\le n+1$ , i.e $f:\mathbb C\to (X,\Delta_{\ell_1})$ is an entire orbifold curve,  then there exists a finite morphism $\pi: X\to \mathbb P^n({\Bbb C})$.  Let $D_0 =\mathbb P^n({\Bbb C})\setminus \mathbb G_m^n$, i.e. $D_0=[x_0=0]+\cdots+ [x_n=0]$.  Let  $R\subset X$  be the ramification divisor of $\pi$ omitting components from the support of $\pi^*(D_0)$, which contains the support of $D$.  Finally, we claim that $\pi(R)$ and $D_0$ are in general position on $\mathbb P^n$.  Since this  is a local property, by the definition of normal crossing, it reduces to a finite morphism $\mathbb A^n \to \mathbb A^n$, which extends to a finite morphism  
 $\mathbb P^n({\Bbb C}) \to \mathbb P^n({\Bbb C})$.  The claim then follows from the arguments in the proof of \cite[Theorem 1.2]{GSW22}.

Since $(X, D)$ is of log-general type implies that we can find a sufficiently large integer $\ell_2\ge \ell_1$ such that $(X,\Delta_{\ell_2})$  is   of  general type.  $($See \cite[Corollary 2.2.24]{LazarsfeldI}.$)$  We can now apply Theorem \ref{finitemorphismtoric} to find a positive integer $\ell\ge \ell_2$ such that if $f:\mathbb C\to (X,\Delta_{\ell})$ is an entire orbifold curve, then $f$ is algebraically degenerate, which gives a contradiction.

 \end{proof}

\end{document}